\documentclass[11pt]{article}

\usepackage{amsmath,amssymb}

\def\tabaddress#1{{\small\it\begin{tabular}[t]{c}#1 \\[1.2ex]\end{tabular}}}
\def\UPCMAT{Departamento de Matem\'atica Aplicada IV\\
Universitat Polit\`ecnica de Catalunya-BarcelonaTech\\
Campus Norte, Ed. C-3. C/ Jordi Girona 1, E-08034 Barcelona, Spain\\
{\rm e-mail:} arturo@ma4.upc.edu}
\def\UPCMA4{Departamento de Matem\'atica Aplicada IV\\
Universitat Polit\`ecnica de Catalunya-BarcelonaTech\\
Campus Norte, Ed. C-3. C/ Jordi Girona 1, E-08034 Barcelona, Spain\\
{\rm e-mail:} matmcml@ma4.upc.edu , nrr@ma4.upc.edu}
\def\UC3M{Department of Mathematics\\
University of  California at Berkeley. Berkeley, CA 94720, USA\\
(On leave of absence from the Departamento de Matem\'aticas\\
Universidad Carlos III de Madrid\\
Avda. de la Universidad 30, 28911 Legan\'es, Madrid, Spain)\\
{\rm email}: albertoi@math.uc3m.es}

\newtheorem{theorem}{Theorem}[section]
\newtheorem{corollary}{Corollary}
\newtheorem{lemma}[theorem]{Lemma}
\newtheorem{proposition}{Proposition}

\newtheorem{definition}[theorem]{Definition}
\newtheorem{remark}{Remark}

\def\beq{\begin{equation}}
\def\eeq{\end{equation}}
\def\bea{\begin{eqnarray}}
\def\eea{\end{eqnarray}}
\def\beann{\begin{eqnarray*}}
\def\eeann{\end{eqnarray*}}
\def\beasn{\begin{sneqnarray}}
\def\eeasn{\end{sneqnarray}}
\def\ben{\begin{enumerate}}
\def\een{\end{enumerate}}
\def\bit{\begin{itemize}}
\def\eit{\end{itemize}}
\def\proof{( {\sl Proof} )\quad}
\def\qed{\ifvmode\removelastskip\fi
{\unskip\nobreak\hfil\penalty50\hbox{}\nobreak\hfil
\hbox{\vrule height1.2ex width1.2ex}\parfillskip=0pt
\finalhyphendemerits=0 \par\smallskip}}
\def\derpar#1#2{\frac{\partial{#1}}{\partial{#2}}}

\def\feble#1{\mathrel{\mathop =\limits_{#1}}}

\def\moment#1#2#3{{#1}_{#2}, \ldots, {#1}_{#3}}
\def\set#1{\{\,#1\,\}}
 1
\def\vf{{\mathfrak{X}}}
\def\df{{\Omega}}

\def\d{{\rm d}}

\def\Real{\mathbb{R}}

\def\Tan{{\rm T}}
\def\Lie{{\rm L}}
\def\inn{{\bf i}}
\def\Cinfty{{\rm C}^\infty}
\def\supp{{\rm supp}}

\def\ds{\displaystyle}

\title{INVARIANT FORMS AND AUTOMORPHISMS OF
 LOCALLY HOMOGENEOUS MULTISYMPLECTIC MANIFOLDS}
\author{
\centerline{\scshape Arturo Echeverr\'\i a-Enr\'\i quez}
\\ 
 \tabaddress{\UPCMAT}
 \\
\centerline{\scshape Alberto Ibort}
\\
 \tabaddress{\UC3M}
\\
\centerline{\scshape Miguel C. Mu\~noz-Lecanda and Narciso Rom\'an-Roy}
\\
   \tabaddress{\UPCMA4}
}
\date{[{\sl J. Geom. Mech.} {\bf 4}(4) (2012) 397--419. doi:10.3934/jgm.2012.4.397 ]}

\bigskip

\begin{document}
\maketitle
\thispagestyle{empty}

\begin{abstract}
It is shown that the geometry of locally homogeneous multisymplectic manifolds
(that is, smooth manifolds equipped with a closed nondegenerate form of degree
$> 1$, which is locally homogeneous of degree $k$ with respect to
a local Euler field) is characterized by their automorphisms.
 Thus, locally homogeneous multisymplectic manifolds extend the
family of classical geometries possessing a similar property: symplectic,
volume and contact.  The proof of the first
 result relies on the characterization
 of invariant differential forms with respect to the graded Lie algebra of
infinitesimal automorphisms, and on the study of the local properties of
Hamiltonian vector fields on locally multisymplectic manifolds.
In particular it is proved that the
group of multisymplectic diffeomorphisms acts (strongly locally) transitively on the manifold.
 It is also shown that the graded Lie algebra of infinitesimal automorphisms
 of a locally homogeneous multisymplectic manifold characterizes
 their multisymplectic diffeomorphisms.
\end{abstract}

\bigskip

\begin{center}
\vbox{AMS s.\,c.\,(2000):
 53C15,53D35, 57R50, 57S25, 58A10.}\null
{\it Key words}: {\sl Multisymplectic manifolds, multisymplectic
diffeomorphisms, invariant forms, Hamiltonian (multi-) vector fields, graded
Lie algebras.}
\end{center}

\section{Introduction}

It is well-known that some classical geometrical structures are determined
by their automorphism groups; for instance it was shown by Banyaga
 \cite{Ba86,Ba88,Ba97} that the geometric structures
 defined by a volume or a symplectic
form on a differentiable manifold are determined by their automorphism groups,
 the groups of volume-preserving and symplectic diffeomorphisms respectively,
  i.e., if $(M_i,\alpha_i)$,
$i=1,2$ are two paracompact connected smooth manifolds equipped with volume or
 symplectic forms $\alpha_i$ and $G(M_i,\alpha_i)$ denotes the group of volume
 preserving or symplectic diffeomorphisms, then if
$\Phi\colon G(M_1,\alpha_1) \to G(M_2,\alpha_2)$ is a group isomorphism, there
 exists (modulo an additional condition in the symplectic case) a unique
 $C^\infty$-diffeomorphism
$\varphi \colon M_1 \to M_2$ such that
$\Phi (f)=\varphi \circ f \circ \varphi^{-1}$,  for every $f\in G(M_1,\alpha_1)$
and $\varphi^*\alpha_2 = c \, \alpha_1$, with $c$ a constant.  In other words,
 group isomorphisms of automorphism groups of classical structures
 (symplectic, volume) are inner, in the sense that they correspond
 to conjugation by (conformal) diffeomorphisms.

An immediate consequence of the previous theorem
is that, if $(M,\alpha )$ is a manifold with a classical structure (volume or symplectic),
and a differential form
$\beta$ is an invariant for the group $G(M,\alpha )$
of volume preserving or symplectic diffeomorphisms,
then necessarily $\beta$ has to be a
constant multiple of exterior powers of $\alpha$.
In other words, the only differential invariants
of the groups of classical diffeomorphism are multiples of exterior powers
of the defining geometrical structure.
The infinitesimal counterpart of this result was already known
in the realm of Classical Mechanics.
 In 1947 Lee Hwa Chung stated a theorem  concerning
the uniqueness of invariant integral forms
 (the Poincar\'e-Cartan integral invariants)
under canonical transformations  \cite{Hw-47}.
 His aim was to use that result in order to
characterize canonical transformations in the
 Hamiltonian formalism of Mechanics;
that is, canonical transformations are characterized as those
transformations mapping every Hamiltonian system into another Hamiltonian
 one with respect to the same symplectic structure.  Afterwards,
 this result was discussed geometrically \cite{LlR-88} and generalized to
presymplectic Hamiltonian systems \cite{CGIR-85,GLR-84}.
 The main result there was that in a given a symplectic (resp. presymplectic)  manifold,
the only differential forms invariant with  respect to all
 Hamiltonian vector fields are multiples of (exterior powers of)
the symplectic (resp. presymplectic) form.
  Since symplectic and presymplectic manifolds represent
 the phase space of regular and singular
Hamiltonian systems, respectively, this result allows one to identify canonical
transformations in the Hamiltonian formalism of Mechanics with the
symplectomorphisms and presymplectomorphisms group, in each case.

Returning to the general problem of the relation between
 geometric structures and
their group of automorphisms, it is an open question
to determine which geometrical
structures are characterized by these groups.  Apart from symplectic and volume,
contact structures also fall into this class \cite{Ba95}.
Grabowski proved that similar statements hold for Jacobi and Poisson manifolds
\cite{At90,Gr93}.  Our main result shows that
 locally homogeneous multisymplectic manifolds are
determined by their automorphisms (finite and infinitesimal).

 Multisymplectic manifolds are one of the natural generalizations
of symplectic manifolds.
A multisymplectic manifold of degree $k$ is a smooth manifold
$M$ equipped with a closed nondegenerate form $\Omega$
 of degree $k\geq 2$ (see \cite{Ca96a,Ca96b} for more  details
on multisymplectic manifolds).
In particular, multisymplectic manifolds include symplectic and volume manifolds.
A diffeomorphism $\varphi$ between two
multisymplectic manifolds $(M_i, \Omega_i)$, $i=1,2$,
will be called a multisymplectic
diffeomorphism if $\varphi^* \Omega_2 = \Omega_1$.
The group of multisymplectic
diffeomorphisms of a multisymplectic manifold $(M, \Omega )$ will be
denoted by $G(M, \Omega )$.
Multisymplectic structures represent distinguished cohomology
classes of the manifold $M$, but their origin as a
geometrical tool can be traced back to the
foundations of the calculus of variations.
It is well known that the suitable geometric
framework to describe (first-order) field theories
are certain multisymplectic manifolds
(see, for instance \cite{CCI-91,EMR-00,GMS-97,Go-91b,HK-2002,Ka-98,
KT-79, LMS-2004, MS-99, PR-2002,Sd-95,Sa-89} and references quoted
therein). The automorphism groups of multisymplectic manifolds
play a relevant role in the description of the corresponding system, and
it is a relevant problem to  characterize
them in similar terms as in symplectic geometry
(see, for instance, \cite{Sh-2011}).
 However many of the multisymplectic structures that arise in applications
are natural generalizations of the canonical symplectic structure on a cotangent bundle.
Hence they posses a local canonical form, and similar properties to the symplectic case
 would be expected to hold.

In this paper we show that this is actually the case for
an important class of multisymplectic forms; namely, those
that satisfy a local homogeneity property similar to the
local homogeneity of symplectic structures.
In fact it is simple to show that given a symplectic form $\omega$ on a manifold $M$,
 locally there always exists a Euler-like vector field $\Delta$
such that $\Lie (\Delta)  \omega = 2\omega$.   Given a point $x_0\in M$
and an open set $U$ (which we may consider to be contained on a local chart)
a vector field $\Delta$ defined on $U$ is said to be an {\sl Euler-like vector field}
 if there are local coordinates $x^i$ defined on $U$ centered at $x_0$
  such that $\Delta_U = x^i \partial /\partial x^i$.

Then, if $(M,\Omega)$ is a multisymplectic manifold with $\Omega$ of degree $k$,
we will say that $(M,\Omega)$ is locally homogeneous if, for evey $x_0\in M$
and any open neighborhood $U$ there exists an Euler--like vector field 
$\Delta$ on it such that $\Lie(\Delta)\Omega = f \Omega$ (see Section \ref{lhom},
definitions \ref{local_euler} and \ref{homogeneous_forms}).
Notice that the canonical multisymplectic structures on multicotangent bundles
are locally homogeneous and, hence, our study should be of interest, for instance,
in the realm of Hamiltonian field theories, where the extended multimomentum phase space
is a multicotangent bundle \cite{CCI-91,EMR-00,ELMR-2005,MS-99,RRSV-2011}.

Then, our main results are:

\begin{theorem}\label{main}
Let $(M_i,\Omega_i)$, $i=1,2$, be two locally homogeneous
multisymplectic manifolds of degree $k$ and let $G(M_i,\Omega_i)$
denote their corresponding
groups of automorphisms.
Let $\Phi\colon G(M_1,\Omega_1)\to G(M_2,\Omega_2)$
be a group isomorphism which is also a homeomorphism when
$G(M_i,\Omega_i)$ are endowed with the point-open topology.
 Then, there exists a $C^\infty$ diffeomorphism
$\varphi\colon M_1\to M_2$, such that $\Phi (f) = \varphi \circ f \circ
\varphi^{-1}$ for all $f\in G(M_1,\Omega_1)$ and the tangent map $\varphi_*$
maps locally Hamiltonian vector fields of $(M_1,\Omega_1)$
 into locally Hamiltonian
vector fields of $(M_2,\Omega_2)$.  In addition,
 if we assume that $\varphi_*$ maps all
infinitesimal automorphisms of $(M_1, \Omega_1)$
 into infinitesimal automorphisms of
$(M_2, \Omega_2)$ then there is a constant $c$ such that
 $\varphi^*\Omega_2 = c\,\Omega_1$.
 \end{theorem}

This result generalizes the main theorems in \cite{Ba86}
 (Theorems 1 and 2), which are in turn generalizations
 of a theorem by Takens \cite{Ta79}.
Contrary to the proof in \cite{Ba86}, the proof presented here
does not rely on the generalization by
Omori \cite{Om74} of Pursell-Shanks theorem
\cite{Pu54}, which do not apply to this situation because of
 the lack of local normal forms for multisymplectic structures.  However, we
use the following partial generalization of Lee Hwa Chung theorem.

\begin{theorem}\label{invar}
  Let $(M,\Omega)$ be a locally homogeneous multisymplectic
manifold of degree $k$; then the only differential forms  of degree $k$
invariant under the graded Lie algebra of infinitesimal automorphisms of
$\Omega$ are real multiples of $\Omega$.
\end{theorem}

Local properties of multisymplectic diffeomorphisms
of  locally homogeneous multisymplectic manifolds
will play a crucial role throughout the discussion.
They stem from a localization property for Hamiltonian vector
fields that will be discussed in Lemma \ref{localization}.
 These local properties are used to prove
a further result that is interesting in itself:
 the group of multisymplectic diffeomorphisms acts
transitively on the underlying manifold.
This result shows that not all multisymplectic manifolds are locally homogeneous.
In fact, R. Bryant \cite{Br87} showed the existence of a multisymplectic structure of
degree 3 such that its automorphism group is the exceptional group $G_2$ which
is not strongly locally transitive on the underlying manifold
(see also \cite{Ib-2000}).

The paper is organized as follows: in Section 2 we
establish some basic definitions and results, mainly related with
the geometry of multisymplectic manifolds.
Section 3 is devoted to the definition and
basic properties of the graded Lie algebra of
the infinitesimal automorphisms of multisymplectic manifolds.
In Section 4, the definition and some characteristics of
locally homogeneous multisymplectic manifolds is stated;
in particular, the localization lemma
for multisymplectic diffeomorphisms, and the strong local transitivity
of the group of multisymplectic diffeomorphisms
is proved for locally homogeneous multisymplectic manifolds.
In Section 5 we prove the main results on the structure of differential invariants of locally
homogeneous multisymplectic manifolds.
Finally, in Section 6, these results are used to characterize the multisymplectic transformations,
while the proof of the main theorem is completed in Section 7.
The paper ends with an appendix where multisymplectic
manifolds admiting Darboux type coordinates are analised
(according to \cite{LMS-2003}), as examples
of locally homogeneous multisymplectic manifolds.

All manifolds are real, paracompact, connected and $C^\infty$. All
 maps are $C^\infty$. Sum over crossed repeated indices is understood.

\section{Notation and basic definitions}

Let $M$ be an $n$-dimensional differentiable manifold.
 Sections of $\Lambda^m(\Tan M)$ are called
{\sl $m$-vector fields} (or {\sl multivector fields} of degree $m$)
in $M$, and we will denote by $\vf^m (M)$ the set of
 $m$-vector fields in $M$.
Let $\Omega\in\df^k(M)$ be a differentiable $k$-form in $M$ ($k\leq n$).
For every $x\in M$, the form $\Omega_x$ establishes a correspondence
$\hat{\Omega}_m(x)$ between the set of
$m$-vectors $\Lambda^m (\Tan_x M)$ and the $(k-m)$-forms
$\Lambda^{k-m} (\Tan_x^* M)$ as
$$
\begin{array}{cccccc}
\hat\Omega_{m}(x) &\colon&\Lambda^m (\Tan_x M) &
\longrightarrow & \Lambda^{k-m}(\Tan_x^* M) &
\\
& & v & \mapsto & \inn(v)\Omega_x & .
\end{array}
$$
If $v$ is homogeneous,
$v = v_1\wedge\ldots\wedge v_m$, then $\inn (v)\Omega_x = \inn
(v_1\wedge\ldots\wedge v_m)\Omega_x =
\inn (v_1)\ldots\inn (v_m)\Omega_x$.
 Thus, an $m$-vector field $X\in\vf^m(M)$ defines a contraction
$\inn (X)$ of degree $m$ of the algebra of differential forms in $M$.
We denote
$$
\ker^m\Omega=\{ X\in\vf^m(M)\, \vert\, \hat\Omega_m(x)(X)=0;\
\mbox{\rm for every $x\in M$}\}
$$

The $k$-form $\Omega$ is said to be {\sl $m$-nondegenerate}
(for $1\leq m\leq k-1$) if, for every $x\in M$, the subspace
$\ker\,\hat{\Omega}_{m}(x)$ has minimal dimension.
Thus:
\begin{enumerate}
\item
If
$\scriptstyle{\left(\begin{array}{c} n \\ m \end{array}\right) } \leq
\scriptstyle{\left(\begin{array}{c} n \\ k-m \end{array}\right) }$,
 then $\dim\,(\ker\,\hat\Omega_m(x)) = 0$.
\item
If
$\scriptstyle{\left(\begin{array}{c} n \\ m \end{array}\right) } >
\scriptstyle{\left(\begin{array}{c} n \\ k-m \end{array}\right) }$, then
$\dim\,(\ker\,\hat\Omega_m(x))=\scriptstyle{\left(\begin{array}{c}n\\ m
\end{array}\right) } - \scriptstyle{\left(\begin{array}{c} n \\ k-m
\end{array}\right) }$, for every $x\in M$.
\end{enumerate}
The form $\Omega$ is {\sl strongly nondegenerate}
if it is $m$-nondegenerate for every
$m=1,\ldots ,k-1$. Thus, the $m$-nondegeneracy of $\Omega$
 implies that the map
 $\hat{\Omega}_m\colon\Lambda^m(\Tan M) \to \Lambda^{k-m}(\Tan^*M)$
is a bundle monomorphism in the case of item 1, or
 a bundle epimorphism in the case of item 2.
The image of the bundle
$\Lambda^m(\Tan M)$ by $\hat{\Omega}_m$ will be denoted by $E_m$.
Often, if there is no risk of confusion,
 we will omit the subindex $m$ and denote
$\hat{\Omega}_{m}$ simply by $\hat{\Omega}$.

\begin{definition}
Let $M$ be an $n$-dimensional differentiable manifold and $\Omega\in\df^k(M)$.
The couple $(M,\Omega)$ is said to be a
{\sl multisymplectic manifold} if $\Omega$ is closed and $1$-nondegenerate.
The degree $k$ of the form
$\Omega$ will be called {\sl the degree of the multisymplectic manifold}.
\end{definition}

Thus, multisymplectic manifolds of degree $k=2$
are the usual symplectic manifolds, and manifolds
 with a distinguished volume-form are
multisymplectic manifolds of degree its dimension.
Other examples of multisymplectic manifolds are provided by
 compact semisimple Lie groups equipped with the canonical
cohomology 3-class, symplectic 6-dimensional Calabi-Yau manifolds
 with the canonical 3-class, etc \cite{Ib-2000}.
 Notice that there are no multisymplectic manifolds of degrees 1 or $n-1$
because $\ker \Omega$ is nonvanishing in both cases.

 Apart from those already cited, another very important class of
multisymplectic manifolds is the  {\sl multicotangent bundle}
(see also Section \ref{apend}: appendix on special multisymplectic manfiolds):
let $Q$ be a manifold, and $\pi\colon\Lambda^k(\Tan^*Q)\to Q$
 the bundle of $k$-forms in $Q$. This bundle
is endowed with a canonical $k$-form
$\Theta\in\df^k(\Lambda^k(\Tan^*Q))$ defined as follows:
if $\alpha\in\Lambda^k(\Tan^*Q)$,
 and $\moment{U}{1}{k}\in\Tan_\alpha(\Lambda^k(\Tan^*Q))$, then
 $$
\Theta_\alpha(U_1,\ldots ,U_k)=\inn(\pi_*U_k\wedge\ldots\wedge\pi_* U_1)\alpha
 $$
If $(x^i,p_{\moment{i}{1}{k}})$ is a system of natural coordinates
in $W\subset\Lambda^k(\Tan^*Q)$, then
$$
\Theta\mid_W=p_{\moment{i}{1}{k}}\d x^{i_1}\wedge\ldots\wedge\d x^{i_k}
$$
Therefore, $\Omega=-\d\Theta\in\df^{k+1}(\Lambda^k(\Tan^*Q))$
is a $1$-nondegenerate form. Then the couple
$(\Lambda^k(\Tan^*Q),\Omega)$ is a multisymplectic manifold.

 Multisymplectic structures of degree $\geq 3$ are abundant.
 In fact, as shown in \cite{Ma70}, if $M$  is a smooth
manifold of dimension $\geq 7$,
 then the space of multisymplectic structures of degree
$3\leq k \leq n-3$ is residual.
 However, there is no local classification of
multisymplectic forms, in general, not even in the linear case.

Finally we introduce the Schouten-Nijenhuis bracket.
 If $X\in\vf^m(M)$, the graded bracket
$$
\Lie(X)=[\d , \inn (X)]=\d\inn (X)-(-1)^m\inn (X)\d
$$
where, as usual,  $\d$ denotes the exterior differential on $M$, defines a new
derivative of degree $m-1$.  If
 $X\in\vf^i(M)$, and $Y\in\vf^j(M)$,
the graded commutator of $\Lie (X)$ and $\Lie (Y)$ is another
operation of degree $i + j -2$  of the same type, i.e.,
there exists a $(i+j-1)$-vector denoted by $[X,Y]$ such that,
$$
\Lie ([X,Y])=[\Lie(X),\Lie(Y)].
$$
The $r$-bilinear assignement $X,Y \mapsto [X,Y]$ is called the
{\sl Schouten-Nijenhuis bracket} of $X,Y$.
It is a generalization of the Lie bracket for
multivector fields
(see also \cite{Ca96b,Ma-97} for a slighly different definition).

Let $X$, $Y$ and $Z$ be homogeneous multivectors of degrees $i,j,p$
respectively, then the Schouten-Nijenhuis bracket verifies the following
properties characterizing this bracket:
\ben
\item
$[X,Y] = -(-1)^{(i+1)(j+1)} [Y,X]$.
\item
$[X, Y\wedge Z] =
[X,Y]\wedge Z + (-1)^{(i+1)j} Y\wedge [X,Z]$.
\item
$(-1)^{(i+1)(p+1)}[X,[Y ,Z] ] + (-1)^{(j+1)(i+1)}[Y,[Z,X] ]
 +(-1)^{(p+1)(j+1)}[Z,[X ,Y] ] = 0 $.
\een
The exterior algebra of
multivectors has the structure of an odd Poisson algebra, sometimes also called
a Schouten algebra. This allows us to define a structure of
an odd Poisson graded manifold on $M$ whose sheaf of superfunctions is given by
the sheaf of multivector fields and the odd Poisson bracket is the Schouten
bracket.

\section{The graded Lie algebra of infinitesimal automorphisms
 of a multisymplectic manifold}

 From now on, $(M,\Omega)$ will be a multisymplectic manifold
of degree $k$.

A {\sl multisymplectic diffeomorphism}
 is a diffeomorphism $\varphi\colon M\to M$
such that $\varphi^*\Omega = \Omega$. A {\sl locally Hamiltonian vector field}
 on $(M,\Omega)$ is a vector field $X$ whose flow consists of multisymplectic
diffeomorphisms. It is clear that $X$ is a locally Hamiltonian vector field
 if and only if $\Lie (X) \Omega = 0$, or equivalently,
 $\inn(X) \Omega$ is a closed $(k-1)$-form.
 This fact leads to the following generalization:

\begin{definition}
 Let $X\in\vf^m(M)$ ($m\geq 1$).
\ben
\item
$X$ is said to be a
{\sl Hamiltonian $m$-vector field}
if $\inn (X)\Omega$ is an exact $(k-m)$-form; that is,
there exists $\zeta \in\df^{k-m-1}(M)$ such that
\beq
\inn (X)\Omega =\d\zeta
\label{ham}
\eeq
$\zeta$ is defined modulo closed $(k-m-1)$-forms in $M$
(they are denoted by $Z^{k-m-1}(M)$).  The class
$\bar\zeta\in\df^{k-m-1}(M)/Z^{k-m-1}(M)$
defined by $\zeta$ is called the
{\rm Hamiltonian} for $X$, and every element in this class
is a {\rm Hamiltonian form} for $X$.
\item
$X$ is said to be a
{\sl locally Hamiltonian $m$-vector field}
if $\inn (X)\Omega$ is a closed $(k-m)$-form.
In this case, for every point $x\in M$, there is an open
neighborhood $W\subset M$ and $\zeta \in\df^{k-m-1}(W)$ such that
$$
\inn (X)\Omega = \d\zeta \qquad {\rm (on\ W)}
$$
As in the above item, changing $M$ by $W$, we obtain the
{\rm Hamiltonian} $\bar\zeta\in\df^{k-m-1}(W)/Z^{k-m-1}(W)$ for $X$,
and the {\rm local Hamiltonian forms} for $X$.
\een

Conversely, $\zeta \in\df^p(M)$ (resp. $\zeta \in\df^p(W)$)
is said to be a {\sl Hamiltonian $p$-form} (resp.
a {\sl local Hamiltonian $p$-form}) if there exists a $(k-p-1)$-vector field $X\in\vf^{k-p-1}(M)$ (resp. $X\in\vf^{k-p-1}(W)$)
such that \ref{ham} holds (resp. on $W$).
\end{definition}

We denote by ${\mathcal H}^p(M)$ the set
of Hamiltonian $p$-forms in $M$ and
by ${\tilde{\mathcal H}}^p(M)$ the set of
Hamiltonian $p$-forms modulo closed
$p$-forms,  ${\tilde{\mathcal H}}^p(M) = {\mathcal H}^p(M) / Z^p (M)$. The classes in
${\tilde{\mathcal H}}^p(M)$ will be denoted by
$\bar{\zeta}$, which means by that the class containing the Hamiltonian
 $p$-form $\zeta$. Let
${\tilde{\mathcal H}}^*(M) = \oplus_{p\geq 0}{\tilde{\mathcal H}}^p(M)$.

\noindent{\bf Remarks.}
\begin{itemize}
\item
If $m>k$ the previous definitions are void and all $m$-vector
fields are Hamiltonian.
\item
There are no Hamiltonian forms of degree
higher than $k-2$.
\item
Every $m$-vector field $X\in\ker^m\Omega$ is a
Hamiltonian $m$-vector field with Hamiltonian the zero class.
\item
 Locally Hamiltonian $m$-vector fields $X$ of degree $k-1$
define closed 1-forms $\inn(X) \Omega$,
which have locally associated a smooth function $f$
(up to constants) called the {\sl local Hamiltonian function} of $X$.
\end{itemize}

\begin{lemma}
\label{surjective}
Let $\Omega\in\df^k(M)$ be a closed $m$-nondegenerate form.
\ben
\item
For every differentiable form $\zeta\in\df^{m-1}(M)$ such that
$\scriptstyle{\left(\begin{array}{c} n \\ m \end{array}\right) } \leq
\scriptstyle{\left(\begin{array}{c} n \\ k-m \end{array}\right) }$,
there exists a $(k-m)$-locally Hamiltonian multivector
field $X$ possessing it as local Hamiltonian form,
i.e. such that $\inn(X)\Omega =\d\zeta$.
As a consequence, the differentials of Hamiltonian $(m-1)$-forms
of locally Hamiltonian $(k-m)$-vector fields
span locally the $m$-multicotangent bundle of $M$, $\Lambda^m(\Tan^*M)$.
\item
 If $\scriptstyle{\left(\begin{array}{c} n \\ m \end{array}\right) } \leq
\scriptstyle{\left(\begin{array}{c} n \\ k-m \end{array}\right) }$
the family of locally Hamiltonian $(k-m)$-vector
fields span locally the $(k-m)$-multitangent bundle of $M$, which is
$$
\Lambda^{k-m}(\Tan_x M)={\rm span}\set{X_x \mid \Lie (X)\Omega =0  \ , \
X\in\vf^{k-m}(M)}.
$$
\een
\end{lemma}
\proof
\ben
\item
Let $\Omega$ be $m$-nondegenerate of degree $k$. The map
$\hat{\Omega}_{k-m}$ has its rank in the bundle $\Lambda^m(T^*M)$,
but $\hat{\Omega}_{k-m}$
is the dual map of $\hat{\Omega}_m$ (up to, perhaps, a minus sign).
And as $\hat{\Omega}_m$ is a monomorphism, then
 $\pm\hat\Omega_m^*=\hat{\Omega}_{k-m}\colon
\Lambda^{k-m}(\Tan M)\to\Lambda^m(\Tan^* M)$
is onto. Then, for every $\zeta\in\df^{m-1}(M)$,
$\d\zeta$ defines a section of $\Lambda^m(\Tan^* M)$; hence we can
choose a smooth $(k-m)$-vector field $X$ such that
$\hat{\Omega}_{k-m}(x)(X_x)=\d\zeta(x)$, for every $x\in M$.

Taking a family of coordinate functions $x^i$,
the same can be done locally for a family
of $(m-1)$-forms $x^{i_1}\d x^{i_2}\wedge\ldots\wedge\d x^{i_m}$,
showing in this way that the differentials of Hamiltonians $(m-1)$-forms
span locally the $m$-multicotangent bundle of $M$.
\item
For every $X\in\vf^{k-m}(M)$, $\inn(X)\Omega\in\df^m(M)$.
But, taking into account the above item,
for every $x\in M$, there exists a neighborhood $U\subset M$ such that
$\inn(X)\Omega\vert_U=f^i\d\zeta_i$, where $f^i\in \Cinfty (U)$
and $\zeta_i\in\df^{m-1}(M)$ with $\inn(X_i)\Omega\vert_U=\d\zeta_i$
for some locally Hamiltonian $(k-m)$-vector fields $X_i$. Therefore
$X\vert_U=f^iX_i+Z$, with $Z\in\ker^{k-m}\Omega$; that is,
$\inn (Z)\Omega=0$, so $Z$ are also
locally Hamiltonian $(k-m)$-vector fields and the proof is finished.
\een
 \qed

Notice that for $m=1$, if $k\geq 2$, then
 $n= \scriptstyle{\left(\begin{array}{c} n \\
1 \end{array}\right) } \leq \scriptstyle{\left(\begin{array}{c} n \\k-1
\end{array}\right) }$.
 Thus if $\Omega$ is $1$-nondegenerate,
the above Lemma states that the differentials
of Hamiltonian functions of locally Hamiltonian
$(k-1)$-vector fields span locally the cotangent bundle of $M$
 and that, in its turn,
the family of these $(k-1)$-Hamiltonian multivector fields span locally the
$(k-1)$-multitangent bundle of $M$. However the previous Lemma says nothing
about the Hamiltonian vector fields.
 We analize this question in the following Section.

\begin{proposition}
\label{lie_algebra}
\begin{enumerate}
\item
A $m$-multivector field $X\in\vf^m(M)$ is a locally Hamiltonian
$m$-multivector field if and only if $\Lie (X)\Omega =0$.
\item
If $X,Y$ are locally Hamiltonian multivector fields, then $[X,Y]$ is
a Hamiltonian multivector field with Hamiltonian form $\inn(X\wedge Y) \Omega$.
\end{enumerate}
\end{proposition}
\proof
Item 1 is immediate.

For item2, if $X$, $Y$ are multivector fields
 of degrees $l,m$ respectively we have,
$$\Lie ([X,Y]) \Omega = \Lie (X)\Lie (Y) \Omega - (-1)^{l+m} \Lie (Y) \Lie (X)
\Omega = 0 \ .$$
Furthermore, $\inn([X,Y]) \Omega = \Lie (X)\inn(Y) \Omega -
(-1)^{l+m}\inn(Y)\Lie (X)\Omega = \d (\inn(X)\inn(Y) \Omega )$.
\qed

We denote respectively by $\vf^m_h(M)$ and
$\vf^m_{lh}(M)$ the sets of
Hamiltonian and locally Hamiltonian $m$-vector fields in $M$.  It is
clear by the previous proposition that
$\bigoplus_{m\geq 0}\vf^m_{lh}(M)$ is  a
graded Lie subalgebra of the graded Lie algebra of multivector fields.
We say that an $m$-vector field is {\sl characteristic}
if it belongs to $\ker^m\Omega$.
The set of characteristic $m$-vector fields constitutes a graded Lie
subalgebra of $\bigoplus_{m\geq 0}\vf^m_{lh}(M)$.  Moreover,
 the characteristic multivector fields define a graded ideal
of the graded Lie algebra of Hamiltonian multivector
fields.  We denote the corresponding quotient graded Lie algebra by
${\mathcal V}_H^* (M,\Omega)$, and
$$
{\mathcal V}_H^* (M,\Omega) = \bigoplus_{m\geq 0} {\mathcal V}_H^m(M,\Omega)
\,,\qquad
{\mathcal V}_H^m(M,\Omega) = \vf^m_{lh}(M) / \ker^m{\Omega} .
$$
Notice that again if $m>k$, $\ker^m \Omega = \Lambda^m (TM)$, hence
 ${\mathcal V}_H^m(M,\Omega) = 0$ and ${\mathcal V}_H^1(M,\Omega)\\ = \vf_{lh}(M)$.
  Namely,
$ {\mathcal V}_H^* (M,\Omega) = \bigoplus_{m= 0}^k {\mathcal V}_H^m(M,\Omega)$.

\begin{definition}
 The Lie algebra ${\mathcal V}_H^* (M,\Omega)$ is called the {\sl
infinitesimal graded Lie algebra} of $(M,\Omega)$ or the
 {\sl graded Lie algebra of infinitesimal automorphisms} of $(M,\Omega )$.
\end{definition}

We can translate this structure of the graded Lie algebra to the corresponding
Hamiltonian forms in  a similar way to that in symplectic geometry (see
\cite{Ca96a} for more details on this construction).
In fact, we can define a graded Lie bracket on
${\tilde{\mathcal H}}^*(M)=\oplus_{p\geq 0}{\tilde{\mathcal H}}^p(M)$ as follows:

\begin{definition}
Given $\bar\xi\in\tilde{\mathcal H}^p(M)$, $\bar\zeta\in\tilde{\mathcal H}^m(M)$,
let $X_{\xi}\in\vf^{k-p-1}_h(M)$,
$Y_{\zeta}\in\vf^{k-m-1}_h(M)$ be
their corresponding Hamiltonian multivector fields modulo
$\ker{\Omega}_*$, where $\ker{\Omega}_*=\bigoplus_{j=0}^k\ker^j\Omega$.
The bracket of these Hamiltonian classes
(related to the multisymplectic structure $\Omega$) is the
$(p+m-k+2)$-Hamiltonian class $\{\bar\xi ,\bar\zeta\}$ containing the form,
$$
\{\xi ,\zeta\}= \Omega (X_{\xi},Y_{\zeta})=
\inn (Y_{\zeta})\inn (X_{\xi})\Omega=
\inn (Y_{\zeta})\d\xi=
(-1)^{(k-p-1)(k-m-1)}\inn (X_{\xi})\d\zeta \ .
$$
\end{definition}

In the same way as in the symplectic case, the Poisson bracket
is closely related to the Lie bracket. Now we have:

\begin{proposition}
Given $X_{\xi}\in\vf^p_h(M)$, $Y_{\zeta}\in\vf^m_h(M)$
Hamiltonian multivector fields, let
$\bar\xi\in{\tilde{\mathcal H}}^{k-p-1}(M)$,
$\bar\zeta\in{\tilde{\mathcal H}}^{k-m-1}(M)$ be the corresponding Hamiltonian classes.
Then the Schouten-Nijenhuis bracket $[X_{\xi},Y_{\zeta}]$
is a Hamiltonian $(p+m-1)$-vector field whose Hamiltonian
$(k-p-m-2)$-form is $\{\zeta ,\xi\}$; that is,
$$
X_{\{\zeta ,\xi\}}=[X_{\xi},Y_{\zeta}] \ .
$$
\end{proposition}
\proof
By definition,
$\inn (X_{\{\zeta ,\xi\}})\Omega =\d\{\zeta ,\xi\}$.
Furthermore, as a consequence of Proposition \ref{lie_algebra}
$$
\inn ([X_{\xi},Y_{\zeta}])\Omega = \d\inn (X_{\xi})\inn(Y_{\zeta})\Omega =
\d\{\zeta ,\xi\} \ .
$$
Thus
$\inn (X_{\{\zeta ,\xi\}})\Omega =\inn ([X_{\xi},Y_{\zeta}])\Omega$,
and therefore $X_{\{\zeta ,\xi\}} =[X_{\xi},Y_{\zeta}]$.
\qed

As a consequence of which, we have:

\begin{proposition}
$({\tilde{\mathcal H}}^*(M),\{\ ,\ \})$ is a graded Lie algebra whose grading is
defined by $|\bar \eta | = k- p -1$, $\eta$ being a $p$-Hamiltonian form.
\end{proposition}

\begin{remark}
The graded Lie algebra ${\mathcal V}_H^* (M,\Omega)$
 possesses as elements of degree
zero the Lie algebra of locally Hamiltonian vector fields on $(M,\Omega)$,
which is the Lie algebra of the ILH-group  \cite{Om74}
of smooth multisymplectic
diffeomorphisms.  This suggests the possibility of embracing in a single
structure of supergroup both smooth multisymplectic diffeomorphisms and
infinitesimal automorphisms of a multisymplectic manifold $(M,\Omega)$.  This
can certainly be done by extending to the graded setting
 some of the techniques used to deal with
ILH-Lie groups.
\end{remark}

\section{Locally homogeneous multisymplectic manifolds.
 The group of multisymplectic diffeomorphisms}
\label{lhom}

As mentioned earlier,
in general, it is not true that the locally Hamiltonian vector fields
in a multisymplectic manifold span the tangent bundle of this manifold.
However, there is a simple property, already mentioned
in the introduction, that implies it (among other things).

\begin{definition}\label{local_euler}
 Let $M$ be a differentiable manifold.
 Consider $x\in M$ and a compact set $K$
 such that $x\in\stackrel{\circ}{K}$.
A {\sl local Liouville} or {\sl local Euler-like vector field }at $x$ with respect to $K$ is a vector field $\Delta^x$ on $M$
such that $\supp\,\Delta^x:=\overline{\{ y\in M\, | \, \Delta^x(y)\not= 0\}} \subset K$,
and there exists a diffeomorphism $\varphi \colon\overbrace{{\rm supp}\,\Delta^x}^{\circ} \to \Real^n$
such that $\varphi_*\Delta^x = \Delta$,
where $\ds\Delta=x^i\derpar{}{x^i}$ is the standard Liouville or dilation vector field in $\Real^n$.
\end{definition}

Notice that if $\Delta^x$ is a local Euler-like vector field at $x$ with respect to $K$
and $\lambda$ is a bump function around $x$ with support contained in $\stackrel{\circ}{K}$,
then $\lambda \Delta^{x}$ is also a local Euler-like vector field at $x$,
but now with respect to $\supp\, \lambda\subset K$.
Notice also that if $\Delta_U$ is a vector field defined locally on an open set $U$
(in this case $U \subset K$ for some compact $K$)
such that in some local coordinates it has the form $x^i\derpar{}{x^i}$,
then by choosing a bump function around $x\in U$ with compact support contained in $U$,
the product $\Delta_\lambda = \lambda \Delta_U$
defines a local Euler-like vector field at $x$ with respect to the compact set $\supp \lambda$.

\begin{definition}
\label{homogeneous_forms}
 A differential form $\Omega$ is said to be {\sl locally homogeneous} at $x\in M$ if,
for every open set $U$ containing $x$, there exists a local Euler-like vector field $\Delta^x$ at $x$
with respect to a compact set $K\subset U$ such that
\beq
\Lie({\Delta^x}) \Omega = f \Omega \ ;\quad f\in\Cinfty(U) \ .
\label{homogeneous_form}
\eeq

The form $\Omega$ is locally homogeneous if it is locally homogeneous for all $x\in M$.

A couple $(M,\Omega)$, where $M$ is a manifold and
$\Omega\in\df^k(M)$ is locally homogeneous is called a
{\rm locally homogeneous manifold}.
 \end{definition}

It is obvious from the definition of $\Delta^x$ that, out of $\supp\,\Delta^x$,
the function $f$ vanishes.
In many instances it is possible to see that if $\Omega$ is locally homogeneous at $x$,
then restricting to a smaller open subset, the function $f$ can be chosen to be constant, $f=c$.

For instance, multicotangent bundles are locally homogeneous manifolds.
In fact, let $\alpha_0\in\Lambda^k(\Tan^*Q)$, and
consider local natural coordinates $(q^i,p_{\moment{i}{1}{k}})$
in a small neighborhood of $\alpha_0$. We define
$\ds r^2:=\sum_i(q^i)^2+\sum_i(p_{\moment{i}{1}{k}})^2$ ,
and $\lambda\equiv\lambda(r^2)$ a bump fuction.
Then consider the local Liouville vector field
$\Delta_\lambda:=\lambda\Delta$. Therefore,
a straightforward calculation shows that,
for the natural multisymplectic form $\Omega\in\df^{k+1}(\Lambda^k(\Tan^*Q))$,
 $$
\Lie(\Delta_\lambda)\Omega=[(k+1)\lambda+2r^2\lambda']\Omega \ .
 $$
and the function $f=(k+1)\lambda + 2r^2 \lambda'$ vanishes identically outside of $\supp \lambda$,
and is constant equal to $k+1$ in a smaller ball around $x$ contained in $\supp \lambda$.
Notice that symplectic manifolds, as well as manifolds endowed with volume forms,
are locally homogeneous.

At this point, we can show that Hamiltonian vector fields
in locally homogeneous multisymplectic manifolds can be localized.
This property plays a crucial role in the discussion to follow.

\begin{lemma}\label{localization0}
Let $(M, \Omega )$ be a locally homogeneous manifold. Then, if $\Delta$ is a local Euler-like vector field
for $\Omega$  (equation \ref{homogeneous_form}), then
its flow leaves invariant the subbundle $E_1=Im\,\hat\Omega_1$.
\end{lemma}
\proof
Let $x$ be a point in $M$.  As $M$ is locally contractible,
let $U$ be a small enough open set and $\Delta^x$ an Euler-like vector field at $x$ with respect to
a compact set $K$ contained on it.
(We can assume that $U$ is a coordinate chart with coordinates $x^i$
 centered at $x$ and adapted to $\Delta^x$ as in Definition \ref{local_euler}).
Let $\varphi_s$ be the local flow defined by $\Delta^x$.  Notice that the flow $\varphi_s$ is given by
$\varphi_s (x) = e^s x$ in the previous local coordinate chart.

The local homogeneity property of $\Omega$ with respect to $\Delta$ implies
 that $\varphi_s^*\Omega= f_s \, \Omega$, for some function $f_s$.
 Notice that
$\varphi_s^*\Omega = f_s\, \Omega$ $\Longleftrightarrow$ $\Lie(\Delta) \Omega = \dot{f}_0\, \Omega$,
 with $\dot{f}_s = df_s/ds$.  Actually,
 if  $\Lie(\Delta) \Omega = f\,  \Omega$, then $f_s(x) =
\int_0^s f(\varphi_{\tau}(x)) d\tau$ and $\dot{f}_0 = f$,
 with $\varphi_s$ the local flow of $\Delta$ such that $\varphi_0(x) = x$.

As a consequence, the flow $\varphi_s$ leaves invariant the
subbundle $E_1$. In fact; if $\eta = \inn(X) \Omega \in E_1$,
then $\varphi_s^*\eta \in E_1$ because
$$
\varphi_s^* \eta = \inn((\varphi_{-s})_*X)\varphi_s^*\Omega =
\inn((\varphi_{-s})_*X)(f_{s}\,\Omega)=
\inn(f_{s}(\varphi_{-s})_*X)\Omega \ .
$$

Moreover, for every $\sigma\colon\mathbb{R}\to\mathbb{R}$,
 with appropriate domain and range, the one-parameter
family of local diffeomorphisms $\tilde{\varphi}_s = \varphi_{\sigma(s)}$ will preserve
$E_1$ too.
\qed

\begin{lemma}\label{localization}
Let $X$ be a locally Hamiltonian vector field on a locally homogeneous multisymplectic manifold
$(M,\Omega)$.
Let $x_0$ be a point in $M$, then for each open set $U$ containing $x_0$,
there exists an open neighborhood $V$ of $x_0$ such that $V\subset \bar{V} \subset U$,
with $\bar{V}$ compact, and a locally Hamiltonian vector field $X'$
such that $X'$ coincides with $X$ in $V$ and vanishes identically outside of $U$.
\end{lemma}
\proof
If $X$ is a locally Hamiltonian vector field, then
$\inn(X) \Omega = \eta $, with $\eta$ a closed $(k-1)$-form.

We prove the Lemma in two steps:

a) The closed form
\beq
\inn(X) \Omega = \eta ,
\label{poincare}
\eeq
has a special exact form in the coordinate chart $(U, x^i)$.

Because $\Omega$ is locally homogeneous, given $x_0$ and $U$ there exists
a local Euler-like vector field at $x_0$ with respect to a compact set $K$ contained in $U$.
We will denote by $\Delta$ such vector field and by $\varphi_t$ its flow.  Notice that the vector
field $\Delta$ is complete because its support is compact.

Let us consider the smooth isotopy $\rho_t = \varphi_{\ln t}$ .
The one-parameter family of local maps
$\rho_t$ for  $t\in[0,1]$ define a
strong deformation retraction from $\stackrel{\circ}{K}$ to $x_0$, i.e., $\rho_1 = \varphi_0 =
Id|_{\stackrel{\circ}{K}}$, and
$\rho_0 = \varphi_{-\infty}$ maps $\stackrel{\circ}{K}$ onto $x_0$.
Notice that $ \varphi_{-\infty} (x) = \lim_{t\to \infty} \varphi_t (x) = x_0$ is
well defined because $\Delta$ being a Euler-like vector field with respect to
$x_0$ has a unique fixed point $x_0$ in the interior of $K$ which is just its stable manifold.

Now we can obtain the local Poincar\'e representation of the closed form $\eta$
in equation \ref{poincare}.  Let $\Delta_t$ be the time-dependent vector
field whose
flow is given by $\rho_t$, i.e.,
 \beq\label{moser} \frac{\d}{\d t} \rho_t
= \Delta_t
\circ \rho_t \ .
\eeq
Then,
$$
\frac{\d}{\d t}(\rho_t^* \eta ) = \rho^*_t (\Lie(\Delta_t)\eta\ .
$$
Hence,
$$ \rho_1^* \eta -\rho_0^* \eta = \int_0^1 \frac{\d}{\d
t}(\rho_t^* \eta
)\,\d t= \int_0^1 \rho^*_t (\Lie (\Delta_t) \eta ) \, \d t = \d\int_0^1
\rho_t^*
(\inn(\Delta_t) \eta) \, \d t \ .
$$
Thus,
$$
\eta = \d \int_0^1 \rho_t^* (\inn(\Delta_t) \eta) \d t \ ,
$$
because $\d\eta=0$ and $ \rho_1^* \eta -\rho_0^* \eta=\eta$.
Hence $\eta = \d\zeta$ with $\zeta = \int_0^1 \rho_t^*(\inn(\Delta_t)\eta)\, \d t $,
 in the open set $\stackrel{\circ}{K}$.

Observe that we have  not used that $\eta\in E_1$.

b) Localization of the Hamiltonian vector field $X$.

We will try to localize the vector field $X$ by using a bump function $\lambda$
centered at $x_0$, i.e., we shall choose $\lambda$ such that the
closure of $V= \supp\, \lambda$ will be a compact set contained in
$\stackrel{\circ}{K}$ and with $x_0$ in its interior.
Unfortunately the vector field $\lambda X$ is not locally Hamiltonian in general,
hence we will
proceed by modifying the Hamiltonian form $\zeta$ of $X$ instead. We define a
new vector
field $\Delta_t^\prime$ by scaling the vector field $\Delta_t$ by
$\lambda$, i.e.,
$$\Delta_t^\prime = \lambda \Delta_t .$$
We denote the flow of $\Delta_t^\prime$ by $\rho_t^\prime$.
Notice that the family of local maps $\rho_t^\prime$ will be obtained from $\rho_t$ by
reparametrizing the parameter $t$. Thus if $\rho_t$ were leaving the subbundle $E_1$ invariant,
the same would be true for $\rho_t^\prime$.

Moreover, we can choose the function $\lambda$ such that $\Delta_t (\lambda ) = r (t)$  and
$$
r (t) = \left\{ \begin{array}{l} 0, ~~ \mbox{if} ~~ 0 \leq t \leq 1/3 \\
r(t), ~~ \mbox{which is a positive function such that}
~~\ds\frac{\d}{d t} r(t) > 0, ~~ \mbox{if}
~~ 1/3 < t < 2/3
\\ 1, ~~ \mbox{if} ~~ 2/3 \leq t \leq 1 \, . \end{array} \right.
$$

Then the flow $(\rho_t^\prime)^*$ leaves invariant the
subbundle $E_1 = \hat{\Omega}_1 (\Tan M)$  and $(\rho_t^\prime)^* \eta \in
E_1$ for all $0\leq t \leq 1$.  Again, repeating the computation leading
to equation \ref{moser}, using the vector field $\Delta_t^\prime$ instead, we get
\beq\label{eta_prime}
(\rho_1^\prime)^* \eta - (\rho_0^\prime)^* \eta = \d \int_0^1
(\rho_t^\prime)^*
(\lambda \inn(\Delta_t) \eta ) \, \d t \ .
\eeq

As in the undeformed situation (with $\lambda =1$), $\rho_1^\prime ={\rm Id}$,
however $\rho_0^\prime$ is not a retraction of $\stackrel{\circ}{K}$ onto $x_0$.
Nevertherless, the $(k-1)$-form $\eta^\prime =  \d \int_0^1 (\rho_t^\prime)^*
(\lambda \inn(\Delta_t)
\eta ) \d t$ is in $E_1$, because both $(\rho_0^\prime)^* \eta$ and
$(\rho_1^\prime)^* \eta$, are in $E_1$. Thus there exists a vector field
$X^\prime$ such that
$$\inn(X^\prime ) \Omega = \eta^\prime  \ .$$
The form $\eta^\prime$ is closed by construction, hence $X^\prime$ is
locally Hamiltonian.

Moreover, if $y$ is a point lying in the interior of the
set $\lambda^{-1}(1) \subset V$ then, $\rho_s^\prime (y) = \rho_s (y)$.
Consequently, from equation \ref{eta_prime}, $\eta^\prime (y) =
\eta (y)$ and $X^\prime (y) = X(y)$.  If, on the contrary, $y$ lies outside
the
compact set $\bar{V}$, we have $\rho_t^\prime (y) = y$, for all
$t$, because $\lambda$ vanishes there,
thus $\Delta_t^\prime$ vanishes and the flow
is the identity. Then $\eta^\prime (y) = 0$ and $X^\prime (y) = 0$.
\qed

A far reaching consequence of the localization lemma is the transitivity of the
group of multisymplectic diffeomorphisms.  We will first prove the following
result:

\begin{lemma}\label{transitivity}
Let $(M,\Omega)$ be a locally homogeneous multisymplectic manifold.
Then the family of locally Hamiltonian vector
fields span locally the tangent bundle of $M$, that is
$$
\Tan_x M = {\rm span} \set{X_x \mid X\in\vf (M), \Lie (X) \Omega = 0} \ .
$$
\end{lemma}
\proof
We will work locally.  Let $U$ be a contractible open neighborhood of a
given point $x\in M$.  We can shrink $U$ to be contained in a coordinate chart
with coordinates $x^i$.  The tensor bundles of $M$ restricted to $U$ are
trivial.  In particular, the subbundle $E_1$ restricted to $U$ is trivial.
Let $v\in \Tan_x M$ be an arbitrary tangent vector.  Let $\nu =
\hat{\Omega}_1(x) v \in E_1 \subset\Lambda^{k-1} (\Tan_x^* M)$.  Consider a
vector field $X$ on $U$ such that $X(x) = v$.  Then $\inn(X) \Omega = \eta$ and
the  $(k-1)$-form $\eta$ is not closed in general.  As in
Lemma \ref{localization} we consider a strong deformation retraction $\rho_t$ and the
corresponding  vector field $\Delta_t$.  Now, first we have
$$ \int_0^1 \frac{\d}{\d s} (\rho_s^* \eta )\,\, \d s = \eta ,$$
and furthermore,
\beq\label{intes}
\int_0^1 \frac{\d}{\d s} (\rho_s^* \eta )\, \d s = \d \int_0^1 \rho_s^*
(\inn(\Delta_s) \eta )\, \d s + \int_0^1 \rho_s^* (\inn(\Delta_s) \d\eta )\, \d s \ .
\eeq
However, as $\d\eta = \d\inn(X) \Omega = \Lie (X) \Omega$, then
$$
\inn(\Delta_t) \d\eta = \inn(\Delta_t)\Lie (X) \Omega =
\inn([\Delta_t, X])\Omega + \Lie (X) \inn(\Delta_t)\Omega \ .
$$
Thus, returning to equation \ref{intes}, we obtain
$$
\eta = \d \int_0^1 \rho_s^* \inn(\Delta_s) \eta\, \d s + \int_0^1 \rho_s^*
(\inn([\Delta_s, X]) \Omega + \Lie (X) \inn(\Delta_s)\Omega ) \, \d s \  .
$$
Choosing the vector field $X$ such that its flow leaves $E_1$ invariant,
the second term on the r.h.s. of the previous equation will be in $E_1$,
so the first term will be in $E_1$ too.  Let us define $\ds\eta' =
\d \left(\int_0^1 \rho_s^* \inn(\Delta_s) \eta\, \d s \right)$, and let us denote by $X'$
the Hamiltonian vector field on $U$ defined by
$$\inn(X') \Omega = \eta'  \ .$$
Evaluating $\eta'$ at $x$ we find that $\eta' (x) = \eta (x)$, hence $X' (x) =v$.
We then localize the vector field $X'$ in such a way that the closure of its support is
compact and is contained in $U$.  We can then extend this vector field
trivially to all $M$, and this extension is locally Hamiltonian.  Finally the
value of this vector field at $x$ is precisely $v$.
\qed

Recall that a group of diffeomorphisms $G$ is said to act
$r$-transitively on $M$ if for any pair of collections $\{ x_1, \ldots, x_r \}$,
$\{ y_1, \ldots, y_r \}$ of distinct points
of $M$, there exists a diffeomorphism $\phi
\in G$ such that $\phi (x_i ) = y_i$.  If the group $G$ acts transitively for
all $r$, then it is said to act $\omega$-transitively or
 transitively for short.
The transitivity of a group of diffeomorphisms
can be reduced to a local problem
because (strong) local transitivity implies transitivity.  More precisely, we
will say that the group of diffeomorphisms $G$ is strongly locally transitive
on $M$ if for each $x\in M$ and a neighborhood $U$ of $x$, there are
neighborhoods $V$ and $W$ of $x$ with
$\bar{V} \subset W \subset \bar{W} \subset
U$, $\bar{W}$ compact, such that for any $y\in V$ there is a smooth
isotopy $\phi_t$ on $G$ joining $\phi$ with the identity, $\phi_1  = \phi$,
$\phi_0 = id$, such that $\phi_1 (x) = y$ and $\phi_t$ leaves every point
outside $\bar{W}$ fixed.  Thus, if $G$ is strongly locally transitive on $M$, then
$G$ acts transitively on $M$ \cite{Bo69}.

\begin{theorem}\label{trans}  The group of multisymplectic diffeomorphisms
$G(M,\Omega )$ of a locally homogeneous multisymplectic manifold
 is strongly locally transitive on $M$.
\end{theorem}
\proof
By Lemma \ref{transitivity} we can construct a local basis of the
tangent bundle in the neighborhood of a given point $x$ made of locally
Hamiltonian vector fields $X_i$.
Using Lemma \ref{localization}, we can localize
the vector fields $X_i$ in such a way that the localized Hamiltonian vector
fields
$X_i^\prime$ will have common supports.
 We denote this common support by  $V$ and
assume that it will be contained in a compact subset contained in $U$.
However, the vector fields $X_i^\prime$ will generate the module
of vector fields inside the support
$V$, so the flows of local vectors fields cover
the same set as the flows of local
Hamiltonian vector fields, although the group of diffeomorphisms is locally strongly
transitive and the same will happen for the group of multisymplectic
diffeomorphisms.
\qed

\begin{corollary}\label{cor_trans}
The group of multisymplectic diffeomorphisms $G(M,\Omega)$ of a
locally homogeneous  multisymplectic
manifold $(M,\Omega)$ acts transitively on $M$.
\end{corollary}
\proof
The conclusion follows from the results in \cite{Bo69} and Theorem \ref{trans}.
\qed

\begin{remark}
The center of the graded Lie algebra $({\tilde{\mathcal H}}^*(M),\{\ ,\ \})$ is a
graded Lie subalgebra, whose elements are called Casimirs.
We must point out that, on locally homogeneous multisymplectic manifolds,
 there are no Casimirs of degree 0, i.e., functions commuting with
anything, because if this was the case,
there would be a function $S$ such that
$ \{ S, \eta \} = 0 $,
for every Hamiltonian form $\eta$.  In particular, $S$ commutes with
$(k-1)$-Hamiltonian forms, but this implies that $X(S) = 0$,
 for every Hamiltonian vector field $X$.  But this is clearly impossible because,
for those kinds of multisymplectic manifolds, Hamiltonian vector fields span
 the tangent bundle by Lemma \ref{transitivity}.
\end{remark}

\section{Invariant differential forms}

In order to prove the main statement in this section,
we first establish two lemmas:

\begin{lemma}
Let $(M,\Omega )$ be a multisymplectic manifold of degree $k$
and $\alpha\in\df^p(M)$ (with $p\geq k-1$) a differential form which
is invariant under the set of
locally Hamiltonian $(k-1)$-vector fields, that is,
$\Lie (X)\alpha =0$, for every $X\in\vf^{k-1}_{lh}(M)$. Then:
\ben
\item
For every $X,Y\in\vf^{k-1}_{lh}(M)$,
\beq
\inn (X)\Omega\wedge\inn (Y)\alpha+\inn (Y)\Omega\wedge\inn (X)\alpha=0
\label{clave1}
\eeq
\item
In particular, for every $X\in\vf^{k-1}_{lh}(M)$ with $\inn (X)\Omega =0$
(that is, $X\in\ker^{k-1}\Omega$), we have
\beq
\inn (X)\alpha =0
\label{clave2}
\eeq
\een
\label{basico1}
\end{lemma}
\proof
\ben
\item
Since $\alpha$ is invariant under $\vf_{lh}^{k-1}(M)$,
for every $X\in\vf_{lh}^{k-1}(M)$, we have $\Lie (X)\alpha  = 0$. Then,
\beq
\d\inn (X)\alpha =(-1)^{k-1}\inn (X)\d\alpha \ .
\label{uno}
\eeq
Given $X,Y\in\vf_{lh}^{k-1}(M)$,
for every $x\in M$ there exists an open neighborhood of $x$,
$U\subset M$, and $f,g\in \Cinfty (U)$ such that
$\inn (X)\Omega\feble{U}\d f$ and $\inn (Y)\Omega\feble{U}\d g$
(from now on we will write $X\vert_U\equiv X_f$ and
$Y\vert_U\equiv X_g$).
Then, consider the locally Hamiltonian vector field
$X_h\in\vf_{lh}^{k-1}(U)$ whose expression in $U$ is
$X_h\feble{U}fX_g+gX_f$; its Hamiltonian function in $U$ is
$h=fg\in\Cinfty (U)$, since
$$
\inn (X_h)\Omega \feble{U}
\inn (fX_g+gX_f)\Omega =f\inn (X_g)\Omega +g\inn (X_f)\Omega =
f\d g+g\d f = \d h \ .
$$
Hence,
$$
\inn (X_h)\alpha\feble{U} f\inn (X_g)\alpha +g\inn (X_f)\alpha \ ,
$$
and then
$$
\d\inn (X_h)\alpha\feble{U}\d f\wedge\inn (X_g)\alpha+f\d\inn (X_g)\alpha+
\d g\wedge\inn (X_f)\alpha+g\d\inn (X_f)\alpha \ .
$$
However, taking into account \ref{uno},
$$
\d\inn (X_h)\alpha\feble{U}
f((-1)^{k-1}\inn (X_g)\d\alpha )+g((-1)^{k-1}\inn (X_f)\d\alpha )=
f\d\inn (X_g)\alpha+g\d\inn (X_f)\alpha \ ,
$$
and comparing both results we conclude
\beq
\d f\wedge\inn (X_g)\alpha+\d g\wedge\inn (X_f)\alpha\feble{U}0 \ ,
\label{dos}
\eeq
which is the local expression of equation \ref{clave1}.
\item
Taking $X\in\ker^{k-1}\Omega$ in \ref{clave1},
the equation \ref{dos} gives $\inn (Y) \Omega \wedge \inn (X)\Omega =0$ for
every $Y$.  But because of Lemma \ref{surjective}, this implies that,
$$
\d g\wedge\inn (X)\alpha\feble{U}0 \ ,
$$
for every $g\in \Cinfty (U)$; hence,
$\inn(X)\alpha =0$, for every $X\in\ker^{k-1}\Omega$.
\een
 \qed

\begin{lemma}
Let $(M,\Omega )$ be a multisymplectic manifold of degree $k$
and $\alpha\in\df^p(M)$, with $p=k-1,k$, a differential form
which is invariant under the set of
locally Hamiltonian $(k-1)$-vector fields. Then:
\ben
\item
If $p=k-1$ then $\alpha =0$.
\item
If $p=k$,
there exists a unique $\alpha'\in\Cinfty(M)$ such that
$$
\inn (X)\alpha =\alpha'\inn (X)\Omega \ ,
$$
for every $X\in\vf^{k-1}_{lh}(M)$.
\een
\label{basico2}
\end{lemma}
\proof
The starting point is the equality \ref{clave1}.
Taking $X=Y\not\in\ker^{k-1}\Omega$
(if $X\in\ker^{k-1}\Omega$, then $\inn (X)\alpha=0$ by hypothesis),
we obtain
\beq
\inn (X)\Omega\wedge\inn (X)\alpha =0 \ ,
\label{tres}
\eeq
for every $X\in\vf^{k-1}_{lh}(M)$. Therefore we have:
\ben
\item
If $p=k-1$, then $\inn (X)\alpha\in \Cinfty (M)$, and
according to the first item of Lemma \ref{surjective}
(for $1$-nondegenerate forms), equation
\ref{tres}, together with \ref{clave2}, leads to the result $\inn (X)\alpha
=0$, for every $X\in\vf_{lh}^{k-1}(M)$.
But, taking into account item 2 of Lemma \ref{surjective}
(for $1$-nondegenerate forms),
this also holds for every $X\in\vf^{k-1}(M)$, and
we must conclude that $\alpha =0$.
\item
If $p=k$ and $\inn (X)\alpha =0$ for all $X$, then $\alpha = 0$.
Thus, let us assume that $\inn (X)\Omega \neq 0$ for some $X$, then the
solution of equation \ref{tres} is
\beq
\inn(X)\alpha =\alpha_X'\inn (X)\Omega \ ,
\label{cuatro}
\eeq
and it is important to point out that
the equation \ref{clave2} for $\alpha$
implies that the function $\alpha'_X\in\Cinfty(M)$ is the same
for every $X,X'\in\vf_{lh}^{k-1}(M)$ such that
$\inn (X)\Omega=\inn (X')\Omega$.

Now, returning to equation \ref{clave1} we obtain the relation
$$
\inn (Y)\Omega\wedge\inn (X)\Omega (\alpha_X'-\alpha_Y')=0 \ .
$$
However, $\alpha_X',\alpha_Y'\in \Cinfty (M)$ are
the unique solution of the respective equations \ref{tres}
for each $X,Y\in\vf_{lh}^{k-1}(M)$;
then we have the following options:
\bit
\item
If $\inn (Y)\Omega\wedge\inn (X)\Omega\not= 0$ then $\alpha_X'=\alpha_Y'$.
\item
If $\inn (Y)\Omega\wedge\inn (X)\Omega=0$ then $X=fY+Z$,
where $f\in\Cinfty (M)$ and $Z\in\ker^{k-1}\Omega$.
Therefore:
\bit
\item
If $X\in\ker^{k-1}\Omega$ then $Y\in\ker^{k-1}\Omega$.
Therefore, taking into account item 2 of Lemma \ref{basico1},
the corresponding equations \ref{cuatro} for $X$ and $Y$
are identities, and thus $\alpha_X'$ and $\alpha_Y'$ are arbitrary functions
which we can take to be equal.
\item
If $X\not\in\ker^{k-1}\Omega$ then $Y\not\in\ker^{k-1}\Omega$.
Therefore, taking into account item 2 of Lemma \ref{basico1},
we have
$$
\inn(X)\alpha=
\inn(fY+Z)\alpha =f\inn(Y)\alpha=
f\alpha_Y'\inn(Y)\Omega =\alpha_Y'\inn (fY+Z)\Omega =\alpha_Y'\inn(X)\Omega \ ,
$$
which, comparing with \ref{cuatro}, gives $\alpha_X'=\alpha_Y'$.
\eit
\eit
In any case, $\alpha_X'=\alpha_Y'$, and
as a consequence, the function $\alpha'$
solution to \ref{tres} is the same for every $X\in\vf_{lh}^{k-1}(M)$.
\een
\qed

At this point we can state and prove the following fundamental result:

\begin{theorem}
\label{pre_invar}
Let $(M,\Omega )$ be a locally homogeneous multisymplectic manifold and
$\alpha\in\df^p(M)$, with $p=k-1,k$, a differential form which
is invariant by the set of locally Hamiltonian
$(k-1)$-vector fields and the set of locally Hamiltonian vector fields;
that is, $\Lie (X)\alpha =0$ and $\Lie (Z)\alpha =0$,
for every $X\in\vf^{k-1}_{lh}(M)$ and $Z\in\vf_{lh}(M)$.
Then we have:
\ben
\item
If $p=k$ then $\alpha =c\, \Omega$, with $c\in\Real$.
\item
If $p=k-1$ then $\alpha =0$.
\een
\end{theorem}
\proof
\ben
\item
First assume that $p=k$.

For every $X\in\vf_{lh}^{k-1}(M)$,
according to Lemma \ref{basico2} (item 2), we have
$$
\inn (X)\alpha=\alpha'\inn (X)\Omega=\inn (X)(\alpha'\Omega) \ ,
$$
where $\alpha'\in\Cinfty (M)$
is the same function for every $X\in\vf_{lh}^{k-1}(M)$.
However, taking into account item 2 of Lemma \ref{surjective}
(for 1-nondegenerate forms), the above equality holds
for every $X\in\vf^{k-1}(M)$, and thus $\alpha=\alpha'\Omega$.

Then, for every $Z\in\vf_{lh}(M)$, by hypothesis
$$
0=\Lie (Z)\alpha = (\Lie (Z)\alpha' ) \Omega .
$$
Hence, $\Lie (Z)\alpha' = 0$, but because of Lemma
\ref{transitivity}, locally Hamiltonian vector fields span the tangent
space, thus $\alpha'=c$ (constant). So
$$
\inn (X)\alpha = c \inn (X)\Omega = \inn (X)(c\Omega ) \ ,
$$
and taking into account item 2 of Lemma \ref{surjective}
(for $1$-nondegenerate forms) again,
this relation also holds for every $X\in\vf^{k-1}(M)$. Hence
we must conclude that
$\alpha =c\Omega$.
\item
If $p=k-1$, the result follows straightforwardly from the first item of
the lemma \ref{basico2}.
\een
 \qed

\begin{remark}

\begin{itemize}
\item
Theorem \ref{invar} is an immediate consequence of Theorem \ref{pre_invar}.
\item
Another immediate consequence of this theorem is that,
if $\alpha\in\df^k(M)$ is a differential form
invariant by the sets of locally Hamiltonian
$(k-1)$-vector fields and locally Hamiltonian vector fields,
then it is also invariant by the set of
locally Hamiltonian $m$-vector fields, for $1<m<k-1$.
\item
As is evident, if $k=2$, we have proved (partially)
the classical {\sl Lee Hwa Chung theorem} for multisymplectic manifolds.
\end{itemize}
\end{remark}

\section{Characterization of multisymplectic transformations}

Now we use the theorems above in order to give several
characterizations of multisymplectic transformations
in the same way as Lee Hwa Chung's theorem
allows us to characterize symplectomorphisms
in the symplectic case \cite{Hw-47,LlR-88,GLR-84}.

A vector field $X$ on a multisymplectic manifold $(M,\Omega )$
 is said to be a {\sl conformal Hamiltonian vector field}
if there exists a function $\sigma$ such that
$$
 \Lie (X) \Omega = \sigma \Omega .
$$
It is immediate to check that,
if $\Omega^r \neq 0$, $r>1$, then $\sigma$ must be constant.
Then:

\begin{definition}
A diffeomorphism $\varphi \colon M_1 \to M_2$
between the multisymplectic manifolds
$(M_i, \Omega_i)$, $i=1,2$, is said to be a
{\sl special conformal multisymplectic diffeomorphism} if there exists $c\in \Real$, such
that $\varphi^*\Omega_2 = c \Omega_1$.  The constant factor $c$ will be called the
{\sl scale} or {\sl valence} of $\varphi$.
 \end{definition}

Therefore we have:

\begin{theorem}\label{cgm}
Let $(M_i,\Omega_i)$, $i=1,2$, be two locally homogeneous multisymplectic manifolds.
A diffeomorphism $\varphi\colon M_1\to M_2$
is a special conformal multisymplectic diffeomorphism if and only if
the differential map $\varphi_*\colon\vf (M_1)\to\vf (M_2)$
induces an isomorphism between the graded Lie algebras ${\mathcal V}_H^* (M_1,\Omega_1)$,
${\mathcal V}_H^* (M_2,\Omega_2)$.  Then
we have that
$$\varphi_*X_{\xi}=\frac{1}{c}X_{\varphi^{*-1}\xi} \ .$$
In addition, if $X_1\in\vf_h^m(M_1)$
is any Hamiltonian (resp. locally Hamiltonian) multivector field
with $\xi_1\in\df^{k-m-1}(M_1)$ a Hamiltonian form for $X_1$
(resp. locally Hamiltonian in some $U_1\subset M_1$),
and $\varphi_*X_1 = X_2\in\vf_{lh}^m(M_2)$
with Hamiltonian form $\xi_2\in\df^{k-m-1}(M_2)$
(resp. locally Hamiltonian in $\varphi (U_1)=U_2\subset M_2$); then
\beq c\xi_1=\varphi^*\xi_2+\eta \ ,
\label{rfh}\eeq
where $\eta\in\df^{k-m-1}(M_1)$ is a closed form.
In other words, $\varphi^*$ induces an isomorphism between
classes of Hamiltonian forms.
\end{theorem}
\proof
Taking into account Proposition \ref{lie_algebra}, we have:

\qquad ($\Longleftarrow$)\quad
For every $X_1\in\vf_h^m(M_1)$ (resp. $X_1\in\vf_{lh}^m(M_1)$)
we have that
$\varphi_ *X_1=X_2\in\vf_h^m(M_2)$ (resp. $X_2\in\vf_{lh}^m(M_2)$).
In any case $\Lie (X_2)\Omega_2=0$, then we obtain
$$
0=\varphi^*\Lie (X_2)\Omega_2=\Lie (\varphi_*^{-1}X_2)\varphi^*\Omega_2=
\Lie (X_1)\varphi^*\Omega_2 \ .
$$
Therefore, by Theorem \ref{pre_invar}, we have that
$\varphi^*\Omega_2=c\Omega_1$.

\qquad ($\Longrightarrow$)\quad
Conversely, for every $X_{\xi_1}\in\vf_h^m(M_1)$ we have
$\inn (X_{\xi_1})\Omega_1-\d\xi_1=0$.
Then, since $\varphi^*\Omega_2=c\Omega_1$,
we obtain
\beann
0&=&\varphi^{*-1}(\inn (X_{\xi_1})\Omega_1-\d\xi_1)=
\inn (\varphi_*X_{\xi_1})\varphi^{*-1}\Omega_1-\varphi^{*-1}\d\xi_1X
\\
&=&\frac{1}{c}\inn (\varphi_*X_{\xi_1})\Omega_2-\d\varphi^{*-1}\xi_1
\quad\Longleftrightarrow\quad
\inn (\varphi_*X_{\xi_1})\Omega_2-
\d\left(\frac{1}{c}\varphi^{*-1}\xi_1\right)=0 \ .
\eeann
So, $\varphi_*X_{\xi_1}=X_{\xi_2}\in\vf_h^m(M_2)$ and its Hamiltonian form
$\xi_2\in\df^{k-m-1}(M_2)$ is related with $\xi_1$ by equation \ref{rfh}.

In an analogous way, using $\varphi^{-1}$, we would prove that
$\varphi_*^{-1}X_2\in\vf_h^m(M_1)$,
for every $X_2\in\vf_h^m(M_2)$.

The proof for locally Hamiltonian multivector fields is obtained
in the same way, working locally on $U_1\subset M_1$
and $U_2=\varphi (U_1)\subset M_2$.
\qed

As a consequence of the previous theorem there is another
characterization of conformal multisymplectomorphisms.

\begin{corollary}
Let $(M_i,\Omega_i)$, $i=1,2$, be two locally homogeneous
multisymplectic manifolds.
A diffeomorphism $\varphi\colon M_1\to M_2$
is a special conformal multisymplectic diffeomorphism if and only if
for every $U_2\subset M_2$ and
for every $\xi_2\in\df^p(U_2)$ and $\zeta_2\in\df^m(U_2)$
(\,$p,m<k-1$\, ), we have
\beq
\varphi^*\{\xi_2,\zeta_2\} =\frac{1}{c}\{\varphi^*\xi_2,\varphi^*\zeta_2\} \ .
\label{pp}
\eeq
\end{corollary}
\proof
Let $X_{\xi_2}\in\vf^{k-p-1}_h(M_2)$ and
$Y_{\zeta_2}\in\vf^{k-m-1}_h(M_2)$
be Hamiltonian multivector fields having $\xi_2$ and $\zeta_2$
as Hamiltonian forms in $U_2$.

\qquad ($\Longrightarrow$)\quad
We have
\beq
\varphi^*\{\xi_2,\zeta_2\} =\varphi^*\inn (Y_{\zeta_2})\d\xi_2=
\inn (\varphi_*^{-1}Y_{\zeta_2})\varphi^*\d\xi_2 \ .
\label{ecu}
\eeq
However, if $\varphi$ is a conformal multisymplectomorphism (of valence $c$),
according to Theorem \ref{cgm},
$\varphi_*^{-1}Y_{\zeta_2}\in\vf_h^m(M_1)$ and
$$
\inn(\varphi_*^{-1}Y_{\zeta_2})\Omega_1=\frac{1}{c}\d\varphi^*\zeta_2 \ .
$$
that is, $\varphi_*^{-1}Y_{\zeta_2} = \frac{1}{c}Y_{\varphi^*\zeta_2}$.
Therefore, because of equation \ref{ecu}, we conclude
$$
\varphi^*\{\xi_2,\zeta_2\} =\inn (\varphi_*^{-1}Y_{\zeta_2})\varphi^*\d\xi_2=
\frac{1}{c}\inn (Y_{\varphi^*\zeta_2})\varphi^*\d\xi_2=
\frac{1}{c}\{\varphi^*\xi_2,\varphi^*\zeta_2\}
$$

\qquad ($\Longleftarrow$)\quad
Assuming that equation \ref{pp} holds and again using the definition
of Poisson bracket, it can be written as
$$
\varphi^*\inn (Y_{\zeta_2})\d\xi_2=
\inn (\varphi_*^{-1}Y_{\zeta_2})\varphi^*\d\xi_2=
\frac{1}{c}\inn (Y_{\varphi^*\zeta_2})\varphi^*\d\xi_2 \ ,
$$
for every $\xi_2$. Hence we conclude that
$\varphi_*^{-1}Y_{\zeta_2}=\frac{1}{c}Y_{\varphi^*\zeta_2}\in\vf_h^m(M_1)$,
for every $Y_{\zeta_2}\in\vf_h^m(M_2)$,
and again because of Theorem \ref{cgm}, $\varphi$
is a special conformal multisymplectomorphism.
\qed

\section{Proof of the main Theorem}
\label{seven}

We will prove now Theorem \ref{main}:

\noindent{\bf Theorem 1.1.\ } {\it
Let $(M_i,\Omega_i)$, $i=1,2$, be two locally homogeneous
multisymplectic manifolds and $G(M_i,\Omega_i)$
 will denote their corresponding
groups of automorphisms.   Let $\Phi \colon G(M_1,\Omega_1) \to G(M_2,
\Omega_2)$ be a group isomorphism which is also a homeomorphism when
$G(M_i,\Omega_i)$ are endowed with the point-open topology.
 Then, there exists a $C^\infty$ diffeomorphism
$\varphi\colon M_1 \to M_2$, such that $\Phi (f) = \varphi \circ f \circ
\varphi^{-1}$ for all $f\in G(M_1,\Omega_1)$ and the tangent map $\varphi_*$
maps locally Hamiltonian vector fields of $(M_1,\Omega_1)$
 into locally Hamiltonian
vector fields of $(M_2,\Omega_2)$.  In addition,
 if we assume that $\varphi_*$ maps all
infinitesimal automorphisms of $(M_1, \Omega_1)$
 into infinitesimal automorphisms of
$(M_2, \Omega_2)$, then there is a constant $c$ such that
 $\varphi^*\Omega_2 = c\,\Omega_1$.}

\proof
Let $M_i$, $i=1,2$, be two locally homogeneous multisymplectic manifolds and
let $\Phi$ be a group isomorphism from $G(M_1,\Omega_1)$ to $G(M_2,\Omega_2)$,
which is in addition a homeomorphism if we endow $G(M_i,\Omega_i)$ with the
point-open topology.   Then, Corollary \ref{cor_trans} implies that
the group $G(M_i,\Omega_i)$ acts transitively on
$M_i$, $i=1,2$. Hence, by the main theorem in \cite{We54} there exists a bijective map
 $\varphi\colon M_1\to M_2$ such that $\Phi (f) = \varphi \circ f \circ
\varphi^{-1}$.  Moreover, the map $\varphi$ is a conformal multisymplectic
diffeomorphism if $\Phi$ verifies the conditions in Theorem \ref{main}, as the
following argument shows.

\bit
\item $\varphi$ is a homeomorphism.

Let ${\mathcal A} (M)$ be the class of fixed subsets of
$G(M,\Omega)$,
i.e.,
$${\mathcal A} (M) = \set{ {\rm Fix}\, (f) \mid f \in G(M,\Omega ) }, ~~~~~  {\rm Fix}\, (f) =
\set{x\in M\, \mid\, f(x) = x} \ .$$
Let ${\mathcal B} (M)$ be the class of complements of elements of ${\mathcal A} (M)$,
which is
$${\mathcal B} (M) = \set{ B = M - A \mid A \in {\mathcal A} (M) } \ .$$
Hence, ${\mathcal B} (M)$ is a class of open subsets of $M$.  If $B\in {\mathcal B}
(M)$ we can construct a multisymplectic diffeomorphism $g$  such that $B$ is
the interior of $\supp (g)$.   In fact, for any point $x\in M$ and a
neighborhood $U$ of
$x$, it follows from Lemma \ref{localization} that there exists $B\in {\mathcal B}
(M)$ such that $x\in B\subset U$.  Thus,
${\mathcal B} (M)$ is a basis for the topology of $M$.  Moreover, if $f\in
G(M_1,\Omega_1)$, then ${\rm Fix}\,  (\varphi \circ f \circ \varphi^{-1} ) = \varphi
({\rm Fix}\, (f))$, and if
$g\in G(M_2, \Omega_2)$, then ${\rm Fix}\, (\varphi^{-1} \circ g \circ \varphi ) =
\varphi^{-1} ({\rm Fix}\, (g))$.   Hence, $\varphi$, $\varphi^{-1}$ take basic open
sets (in ${\mathcal B} (M)$) into basic open sets. Thus, they are both continuous, i.e.,
$\varphi$ is a homeomorphism.
\item
 $\varphi$ is a smooth diffeomorphism.

To prove this we adapt the proof in \cite{Ta79} and \cite{Ba86} to our
setting.  To prove that $\varphi$, and
$\varphi^{-1}$ are $C^\infty$ it is enough to show that $h\circ \varphi \in
C^\infty (M_1)$ for all $h\in C^\infty (M_2)$ and $k\circ \varphi^{-1} \in
C^\infty (M_2)$ for all $k\in C^\infty (M_1)$.

Let $x\in M_1$ and $U$ be an open neighborhood of $x$, which is the domain of a
local coordinate chart $\psi \colon U \to \Real^n$.
According to Lemma \ref{transitivity}, there exist
Hamiltonian vector fields $X_i$, with compact supports on $U$, which
are a local basis for the vector fields on an open neighborhood of $x$
contained in $U$.
Let $\phi_t^i$ be the 1-parameter group of diffeomorphisms
generated by $X_i$.
Now let $X$ be any locally Hamiltonian vector field on $M_1$, which we localize
on a neighborhood of $x$ in such  way that its compact
support will be contained in
$U$.  We will denote the localized vector field again by $X$.  Let $\phi_t$
 be 1-parameter group of multisymplectic
diffeomorphisms generated by $X$ (which
exists because $X$ is complete).  For each $t$, $\Psi_t := \Phi (\phi_t) =
\varphi
\circ \phi_t \circ \varphi^{-1}$ is a $C^\infty$
multisymplectic diffeomorphism. The evaluation map,
$$\begin{array}{ccccl}\Psi &\colon&\Real \times M_2 &
\longrightarrow & M_2 \\ & & (t,x) & \mapsto & \Phi_t (x) = \Phi (\phi_t) (x) =
\varphi \circ \phi_t \circ \varphi^{-1} (x) \end{array} $$
is continuous.  Moreover $\Psi_0 = {\rm Id}$ and
$\Psi_{t+s} = \Psi_t \circ \Psi_s$.  Therefore, the map $\Psi$
is a continuous action of $\Real$ on $M_2$ by $C^\infty$
diffeomorphisms.  By the Montgomery-Zippin theorem,
since $\Real$ is a Lie group, this action is
$C^\infty$, i.e, $\Psi$ is smooth in both variables $t$ and $x$.
  Therefore, the
1-parameter group of multisymplectic diffeomorphisms
 $\Psi_t$ has an infinitesimal
generator, i.e., a $C^\infty$ locally Hamiltonian vector field $X_\Psi$ such that,
$$ \frac{\d}{\d t} \Psi_t = X_\Psi \circ \Psi_t .$$
Given $h\in C^\infty (M_2)$, its directional derivative $X_\Psi (h)$ is a
$C^\infty$ function.  For any $x\in M_1$ we have,
$$ X_\Psi (h) (\varphi (x)) = \left. \frac{\d}{\d t} h(\Psi_t(\varphi
(x))\right|_{t=0} = \left. \frac{\d}{\d t} (h\circ\varphi) (\phi_t (x))
\right|_{t=0} \ .$$
Therefore, if $X$ is any of the Hamiltonian vector fields $X_i$ above,
for all $y$ in a small neighborhood of $x$,
the preceding formula gives
$$
X_i(h\circ\varphi )(y)=
\left. \frac{\d}{\d t} (h\circ \varphi) (\phi_t^i (y)) \right|_{t=0} =
(X_i)_\Psi
(h) (\varphi (y)) \ .
$$
This formula shows that $h\circ \varphi$
is a $C^1$-map and that for every locally
Hamiltonian vector field $X$,
\beq\label{transf} (X_\Psi (h)) \circ \varphi = X (h\circ \varphi ) \ .
\eeq
To compute higher partial derivatives, we just iterate this formula using the
vector fields $X_i$; for instance,
$$
(X_j)_\Psi ((X_i)_\Psi (h))\circ \varphi = X_j (X_i (h\circ \varphi )) \ .
$$
Since the Hamiltonian vector fields $X_i$
are a local basis for the vector fields on an open neighborhood of $x$,
we have proved that $h\circ \varphi \in C^\infty (M_1)$.
\item
  $\varphi_*$ maps locally Hamiltonian vector fields
into locally Hamiltonian vector fields.

Equation \ref{transf} shows that $X_\Psi = \varphi_* X$
and because $\Psi_t$ is a flow of multisymplectic diffeomorphisms, then
$\varphi_*X$ is another locally
Hamiltonian vector field.  Thus,
$\varphi_*$ maps every locally Hamiltonian vector field into a
locally Hamiltonian vector field.
\item
$\varphi$ is a special conformal multisymplectic diffeomorphism.

Finally we show that, with the additional hypothesis stated in Theorem \ref{main}, then
$\varphi^* \Omega_2 = c \Omega_1$.
In fact, if in addition we assume that the tangent map
$\varphi_*$ maps all infinitesimal automorphisms of $(M_1, \Omega_1)$ into
infinitesimal automorphisms of $(M_2, \Omega_2)$, then
as a consequence of Theorem \ref{cgm}, we have that
$\varphi^*\Omega_2 = c\, \Omega_1$.
(It is important to point out that this conclusion cannot be reached
unless this new hypothesis is assumed, since the starting set
of assumptions allows us to prove only that
$\varphi_*$ maps locally Hamiltonian vector fields
into locally Hamiltonian vector fields.
However, this result cannot be extended to
Hamiltonian $m$-multivector fields, with $m>1$.)
\eit
 \qed

\section{Conclusions}

We show that locally homogeneous multisymplectic forms
 are characterized by their
automorphisms (finite and infinitesimal).
  As pointed out in the introduction it is
remarkable that a Darboux-type theorem
 for multisymplectic manifolds is not known,
although a class of multisymplectic manifolds with a local structure defined by
Darboux type coordinates has been characterized \cite{Ca96b}.  This obliges
us to use a proof that does not rely on normal forms.

  The statement in
Theorem \ref{main} can be made slighty more restrictive assuming that we are
given a bijective map $\varphi \colon M_1 \to M_2$
such that it sends elements $f\in G(M_1, \Omega_1)$
to elements $\varphi \circ f \circ\varphi^{-1} \in G(M_2, \Omega_2)$. Then
throughout the
proof of the theorem in Section \ref{seven} we show that $\varphi$ is $C^\infty$.  The
generalization we present here uses the transitivity of the group of
multisymplectic diffeomorphisms and is a simple consequence of theorems by
Wechsler \cite{We54} and Boothby \cite{Bo69}.
However, we do not know yet if the
continuity assumption for $\Phi$ can be dropped and replaced by weaker
conditions as in the symplectic and contact cases.   To answer these
questions, it would be necessary to describe the algebraic structure of the
graded Lie algebra of infinitesimal automorphisms of the geometric structure as in the
symplectic and volume cases \cite{Ba78}.  A  necessary first step in this direction
will be to describe the extension of Calabi's
invariants to the multisymplectic setting.

We wish to stress that in the analysis of multisymplectic
structures beyond the symplectic and volume manifolds, it is necessary to
consider not only vector fields, but the
graded Lie algebra of infinitesimal automorphisms of
arbitrary order. Only the Lie subalgebra of derivations of degree zero is
related to the group of diffeomorphisms, but derivations of
higher degrees are needed to characterize the invariants.

Finally, we wish to remark that Theorem \ref{pre_invar}
(which plays a relevant role in this work) is just a partial
geometric generalization for multisymplectic manifolds
of Lee Hwa Chung's theorem. A complete generalization
would have to characterize invariant forms of every degree.
Our guess is that, in order to achive this, additional hypotheses
must be considered; namely: strong nondegeneracy
of the multisymplectic form and invariance by locally Hamiltonian
multivector fields of every order.
Nevertheless, it is important to point out that the
 hypotheses that we assume here
(1-nondegeneracy, local homogeneity and invariance
 by locally Hamiltonian vector fields and
locally Hamiltonian $(k-1)$-multivector fields) are sufficient
for our aim. This is a relevant fact since, for example,
in the jet bundle description of classical field theories
(the regular case),
the Lagrangian and Hamiltonian multisymplectic forms
are just 1-nondegenerate, and in the analysis of field equations,
only locally Hamiltonian $(k-1)$-multivector fields
are relevant \cite{EMR-97,EMR-99,ELMR-2005,PR-2002}.

\section{Appendix: Locally special multisymplectic manifolds}
\label{apend}

We have seen that multicotangent bundles are locally homogeneous manifolds.
It is remarkable that  a Darboux-type theorem
for multisymplectic manifolds, in general, is not known,
although a class of multisymplectic manifolds with a local structure defined by Darboux type
coordinates was characterized in \cite{Ca96b}.
 In \cite{LMS-2003}, multisymplectic
manifolds admiting Darboux-type coordinates
(or what is equivalent, the existence of normal forms for the multisymplectic structure)
have been classified, and this is a sufficient condition which
guarantees that property. In fact,
this is a large class  of multisymplectic manifolds having the property
of being locally homogeneous.

In order to introduce these manifolds, some previous concepts are required
(see \cite{Ca96b,LMS-2003} for details concerning the
definitions and the subsequent results).
First, given a multicotangent bundle $\Lambda^k(\Tan^*Q)$,
let us denote by $\Lambda^k_r(\Tan^*Q)$ the subbundle
made of those $m$-forms in $Q$ vanishing when applied to
$r$ vector fields in $Q$.
Let $\rho\colon\Lambda^k(\Tan^*Q)\to Q$ and
$\rho_r\colon\Lambda^k_r(\Tan^*Q)\to Q$
be the natural projections, and $\Theta_Q\in\df^k(\Lambda^k(\Tan^*Q))$,
$\Theta^r_Q\in\df^k(\Lambda^k_r(\Tan^*Q))$ the tautological $m$-forms in these bundles.
 Then, following the terminology introduced by Tulczyjew in
 \cite{Tu-76a,Tu-76b}, we define:

\begin{definition}
A {\sl special multisymplectic manifold} is a multisymplectic
manifold $(M,\Omega)$ (of degree $k$) such that:
\ben
\item
$\Omega=\d\Theta$, for some $\Theta\in\df^{k-1}(M)$.
\item
There exist a diffeomorphism $\phi\colon M\to \Lambda^{k-1}(\Tan^*Q)$,
with $\dim\, Q=n\geq k-1$,
(or $\phi\colon M\to \Lambda^{k-1}_r(\Tan^*Q)$),
and a fibration $\pi\colon M\to Q$
such that $\rho\circ\phi=\pi$
(resp. $\rho_r\circ\phi=\pi$),
and $\phi^*\Theta_Q=\Theta$ (resp. $\phi^*\Theta_Q^r=\Theta$).
\een
($(M,\Omega)$ is said to be {\rm multisymplectomorphic} to a bundle of forms).
\end{definition}

It is important to point out that every special multisymplectic manifold
has a local chart of {\sl Darboux coordinates} around every point;
that is, we have coordinates $(x^i;p_{i_1\ldots i_{k-1}})_{i=1...n}$
such that
\beq
\Theta=p_{i_1\ldots i_{k-1}}\d x^{i_1}\wedge\ldots\wedge\d x^{i_{k-1}}
\quad ,\quad
\Omega=\d p_{i_1\ldots i_{k-1}}\wedge\d x^{i_1}\wedge\ldots\wedge\d x^{i_{k-1}} \ .
\label{locomega}
\eeq

If $(M,\Omega)$ is a multisymplectic manifold of degree $k$
and ${\mathcal W}$ a distribution in $(M,\Omega)$,
for every $x\in M$ and $l<k$, we define the vector space
$$
{\mathcal W}(x)^{\perp,l}=\{ v\in\Tan_xM\,\vert\,
\inn(v\wedge w_1\wedge\ldots\wedge w_l)\Omega_x=0,\
\mbox{\rm for every $w_1,\ldots,w_l\in{\mathcal W}(x)$}\} \ .
$$
Then, ${\mathcal W}$ is said to be an {\sl $l$-isotropic distribution} if
${\mathcal W}(x)\subset{\mathcal W}(x)^{\perp,l}$,
for every $x\in M$. Therefore:

\begin{definition}
Let $(M,\Omega)$ be a multisymplectic manifold of degree $k$,
and ${\mathcal W}$ a $1$-isotropic involutive distribution of $(M,\Omega)$.
\ben
\item
The triple $(M,\Omega,{\mathcal W})$ is a
{\sl multisymplectic manifold of type $(k,0)$} if,
for every $x\in M$, we have that:
\ben
\item
${\rm dim}\,{\mathcal W}(x)={\rm dim}\,\Lambda^{k-1}(\Tan_xM/{\mathcal W}(x))^*$.
\item
${\rm dim}\,(\Tan_xM/{\mathcal W}(x))>k-1$.
\een
\item
A {\sl multisymplectic manifold of type $(k,r)$}
($1\leq r\leq k-1$) is a quadruple
$(M,\Omega,{\mathcal W},{\mathcal E})$ such that
${\mathcal E}$ is a ``generalized distribution'' on $M$
(in the sense that, for every $x\in M$,
${\mathcal E}(x)\subset\Tan_xM/{\mathcal W}(x)$ is a vector subspace) and,
for every $x\in M$, denoting by
$\pi_x\colon\Tan_xM\to\Tan_xM/{\mathcal W}(x)$
the canonical projection, we have that:
\ben
\item
$\inn(v_1\wedge\ldots\wedge v_r)\Omega_x=0$, for every $v_i\in\Tan_xM$
such that $\pi_x(v_i)\in{\mathcal E}(x)$ ($i=1,\ldots,r$).
\item
${\rm dim}\,{\mathcal W}(x)={\rm dim}\,\Lambda_r^{k-1}(\Tan_xM/{\mathcal W}(x))^*$,
where the horizontal forms are considered with respect to the
subspace ${\mathcal E}(x)$.
\item
${\rm dim}\,(\Tan_xM/{\mathcal W}(x))>k-1$.
\een
\een
\end{definition}

And the fundamental result is:

\begin{proposition}
Every multisymplectic manifold of type $(k,0)$ (resp. of type $(k,r)$)
is locally multisymplectomorphic to a bundle of $(k-1)$-forms
$\Lambda^{k-1}(\Tan^*Q)$ (resp. $\Lambda^{k-1}_r(\Tan^*Q)$),
for some manifold $Q$; that is, to a canonical multisymplectic manifold.
Therefore, there is a local chart of Darboux coordinates around every point.
\end{proposition}

\begin{definition}
Multisymplectic manifolds which are locally multisymplectomorphic
to bundles of forms are called {\sl locally special multisymplectic manifolds}.
\end{definition}

Obviously, every special multisymplectic manifold is a
locally special multisymplectic manifold.

\begin{remark}
As is evident, locally special multisymplectic manifolds have local
Euler-like vector fields; in particular, the local vector fields
$\ds x^i\derpar{}{x^i}+p_{i_1\ldots i_k}\derpar{}{p_{i_1\ldots i_k}}$.
Thus, the corresponding multisymplectic forms are locally
homogeneous.
\end{remark}

A far reaching consequence of all of this is the transitivity of the
group of multisymplectic diffeomorphisms since,
if $(M,\Omega)$ is a locally special multisymplectic manifold,
then the family of locally Hamiltonian vector
fields span locally the tangent bundle of $M$.
In fact, taking into account the local expression \ref{locomega}
of $\Omega$, we have that the local vector fields
\(\displaystyle \left\{\derpar{}{x^i},\derpar{}{p_{i_1\ldots i_k}}\right\}\) are
locally Hamiltonian.

\section*{Acknowledgments}

The authors wish to thank F. Cantrijn and G. Bor for their useful
remarks and comments during the preparation of this work.
The authors would also thank an anonymous referee for pointing
out an error in a previous draft of this work.

We acknowledge the financial support of the
\emph{ Ministerio de Econom\'\i a y Competitividad},
project  MTM2011-15725-E.
The author AI wishes to ackwnowledge the partial
financial support provided by MTM2010-21186-C02-02 and QUITEMAD programe,
and {\sl Fundaci\'on Caja Madrid}, Programa de Movilidad de Profesores de Universidades P\'ublicas de Madrid.
The authors AEE, MCML, NRR wish to thank the financial support
of the  {\sl Ministerio de Ciencia e Innovaci\'on}, projects
MTM2008-00689, MTM2011-22585,
and AGAUR, project 2009 SGR:1338.

We also thank Mr. Jeff Palmer for his assistance in preparing the English
version of the manuscript.

\end{document}